\pdfoutput=1 
\documentclass[12pt]{article}
\usepackage{amsmath,amsfonts,amssymb,url,color,epsfig,graphics}
\usepackage[latin1]{inputenc}
\usepackage{hyperref}

\newcommand{\R}{{\mathbb{R}}}

\newcommand{\Z}{{\mathbb{Z}}}
\newcommand{\N}{{\mathbb{N}}}

\newcommand{\D}{{\mathbb{D}}}
\newcommand{\T}{{\mathbb{T}}}

\def\ha{\frac{1}{2}}

\def\pa{\partial}
\def\ra{\rightarrow}
\def\ga{\alpha}

\def\ge{\varepsilon}

\def\gl{\lambda}

\newtheorem{prop}{Proposition}[section]

\newtheorem{theo}{Theorem}[section]

\begin{document}

\title{Large-time  asymptotics  of the wave fronts length I \\
 The Euclidean disk}

\author{Yves  Colin de Verdi\`ere\footnote{Universit\'e Grenoble-Alpes,
Institut Fourier,
 Unit{\'e} mixte
 de recherche CNRS-UGA 5582,
 BP 74, 38402-Saint Martin d'H\`eres Cedex (France);
{\color{blue} {\tt yves.colin-de-verdiere@univ-grenoble-alpes.fr}}}\\
 David Vicente, \footnote{Universit\'e d'Orl\'eans, Institut Denis Poisson, Unit{\'e} mixte
 de recherche CNRS 7013, Rue de Chartres,
B.P. 6759,
45067 Orl\'eans Cedex 2 (France)
{\color{blue} {\tt david.vicente154@gmail.com}}}}

\maketitle
In the paper \cite{Vi-20}, the second author proves that the length $|S_t|$ of the wave front $S_t$ at time $t$ of a wave propagating in an
Euclidean  disk $\D $ of radius $1$, starting from a source   $q$, admits a linear asymptotics as $t\ra + \infty $:
$|S_t|=\lambda (q)t + \mathrm{o}(t)$ with $ \lambda (q)=2 \arcsin a$ and $a=d(0,q)$. We will give a more direct proof and compute the oscillating corrections to  this linear asymptotics. The proof is based on the ``stationary phase'' approximation.

\section{Wave fronts}
Let us consider a 2D-Riemannian compact  manifold $(X,g)$ possibly with a smooth convex boundary. 
We denote by $g^\star : T^\star X \ra \R $  the half of the dual metric which is the Hamiltonian of the geodesic flow.

 We denote by $\pi _X$ the canonical projection
of $T^\star X$ onto $X$ and $\phi_t:T^\star X \ra T^\star X,~t\in \R$ the Hamiltonian flow of $g^\star $ which is the geodesic flow.
If $X $ has a non empty boundary, we define $\phi_t $  using the law of reflection.       
Let $q\in X$ be given. For any $t>0$, we define the wave front $S_t$ at time  $t$  as the set of points of $X$ 
of the form $\pi_X (\phi_t (\Sigma ^q ) )$ where  
$\Sigma^q:=\{ (q,\xi)\in T^\star X |g^\star (q,\xi)=1 \}$. The set $S_t$ could also be defined as the image by the exponential map
at $q$ of the circle  $\Sigma _t$ of radius $t$ in the tangent space $T_q X$. 

Let us define the length of $S_t$ and denote it by $|S_t|$.
The wave front $S_t$ is a curve parametrized by a circle: $S_t:={\rm exp}_q (\Sigma _t  )$.
This allows to define its length using the Riemannian metric. Note that $S_t$ can admit some singular points. The length
of the corresponding part vanishes and the remaining part is an immersed co-oriented curve with only transversal self-intersections. 

In this article, we focus on the case where $X$ is the unit disk in $\mathbb{R}^2$ and $g$ is the Euclidean metric. In this context, we will prove that the following expansion holds:
\begin{equation}\label{resultat_principal}
|S_t|= 2\ga _0 t + t\sum _{n=0}^\infty  J_n^{\rm approx}(t)+\mathrm{O}(1)
\end{equation}
as $t\to +\infty$, with
\[ J_n^{\rm approx}(t)= \frac{-8\sqrt{2}}{\pi^2 (2n+1)^{5/2}\sqrt{t}}\cos \left( (2n+1)\pi a \right) \cos \left( 
\pi \left( (2n+1)t +\frac{1}{4}\right) \right)\]
where $a$ is the distance from the point $q$ to the center of the disk.

The case of closed surfaces with integrable geodesic flows will be the subject of \cite{CV-20?}.

\section{Numerics}
\begin{figure}[h!]
\centering
\includegraphics[height=2.9cm]{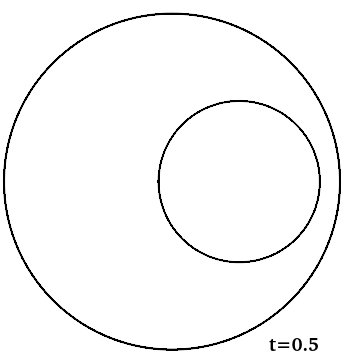} \hspace{0.1cm}
\includegraphics[height=2.9cm]{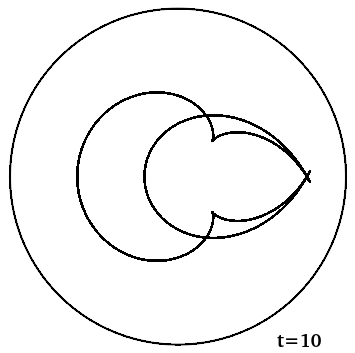} \hspace{0.1cm}
\includegraphics[height=2.9cm]{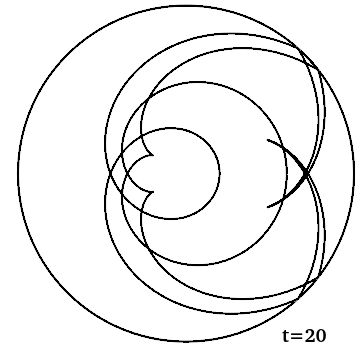} \hspace{0.1cm}
\includegraphics[height=2.9cm]{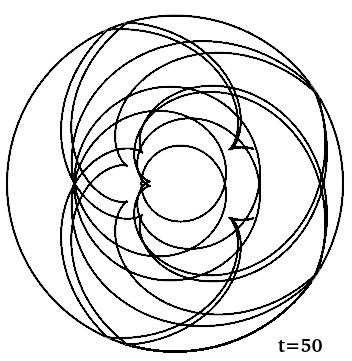}
\caption{Wave Front for $a=0.5$ and $t\in \{0.5,10,20,50\}$}{\label{frontdonde1}}
\end{figure}

In this section, we will compare the expansion given by $(\ref{resultat_principal})$ with the numerical calculations. We introduce a (small) time step $\delta_t>0$, a (large) number of points $n$ which compose the wave front, two vectors $M$ and $V$ in $(\mathbb{R}^2)^n$ such that, for any $k\in[\![1,n]\!]$, $X_k\in\mathbb{R}^2$ represents the position and $V_k\in\mathbb{R}^2$ the speed of the $k$th point of the wave front at a given time. We fix $a\in ]0,1[$ such that $(a,0)$ are the coordinates of the source $q$. Thus, we introduce the following iterative scheme
\begin{equation*}{\label{NUM}}
\left\{
\begin{array}{l}
\mbox{\textbf{initialization:}} \\
\begin{array}{l}
  M \leftarrow \left((a,0), \ldots , (a,0)\right) \in (\mathbb{R}^2)^n,\\
  \mbox{for any } k\in [\![1,n]\!], V_k \leftarrow \left(\cos\left(\frac{2k\pi}{n}\right), \sin\left(\frac{2k\pi}{n}\right)\right),\\
\end{array} \\
\mbox{\textbf{iterative step:}} \\
\begin{array}{l}
  \widetilde{M} \leftarrow M+\delta_t V, \\
  \mbox{for any $k\in [\![1,n]\!]$, }
  \left\{
  \begin{array}{l}
    \mbox{ if } {\widetilde{M}}_k \in \mathbb{D} \mbox{ then } M_k \leftarrow {\widetilde{M}}_k,\\
    \mbox{ else } \\
    \mbox{ } \quad \mbox{compute } {\delta_t^k} \mbox{ s.t. } \Vert M_k+ {\delta_t^k}V_k\Vert = 1 \mbox{ and } {\delta_t^k}\geq 0,\\
    \mbox{ } \quad M_k \leftarrow M_k + \delta_t^kV_k,\\
 \mbox{ } \quad V_k \leftarrow V_k - 2\left< V_k | \frac{M_k}{\Vert M_k\Vert }\right> \frac{M_k}{\Vert M_k \Vert},\\
 \mbox{ } \quad M_k \leftarrow M_k+(\delta_t-\delta_t^k)V_k.
  \end{array}
  \right.\\
\end{array}
\end{array}
\right.
\end{equation*}
The iterative loop consists in the computation of a linear motion outside the boundary and at the boundary one applies the familiar law \textit{the angle of incidence equals the angle of reflection}. After $p$ iterations, $M$ represents the points of the wave front (see Figure \ref{frontdonde1} and \href{https://www.youtube.com/channel/UCMTvpxuhYwbYBYDErSlU0EA/}{\textit{Videos}\footnote{{\color{blue} {\tt https://www.youtube.com/channel/UCMTvpxuhYwbYBYDErSlU0EA/}}}}).

\begin{figure}[!ht]
\centering
\includegraphics[height=5cm]{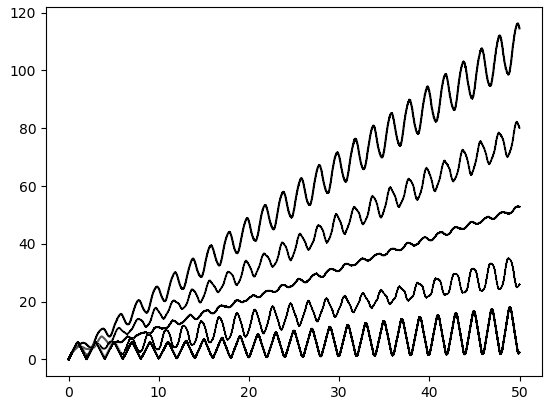}
\caption{Length of the wave front $\vert S_t\vert$ for $t\in [0,50]$ and for differents values of $a$, respectively (from bottom to top), for $a=0.1$, $a=0.3$, $a=0.5$, $a=0.7$ and $a=0.9$}{\label{frontdonde2}}
\end{figure}

\begin{figure}[!ht]{\label{frontdonde3}}
\centering
\includegraphics[height=5cm]{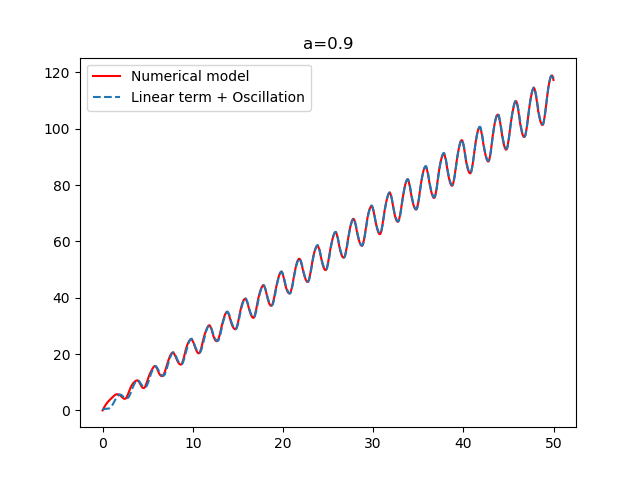} \hspace{0.1cm}
\includegraphics[height=5cm]{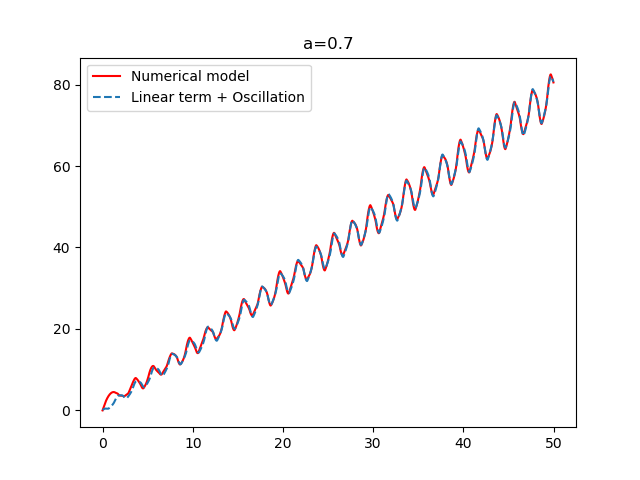} \hspace{0.1cm}\\
\includegraphics[height=5cm]{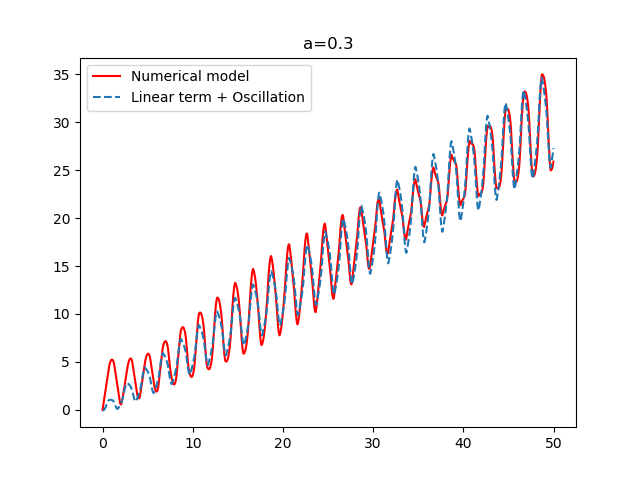}
\hspace{0.1cm}
\includegraphics[height=5cm]{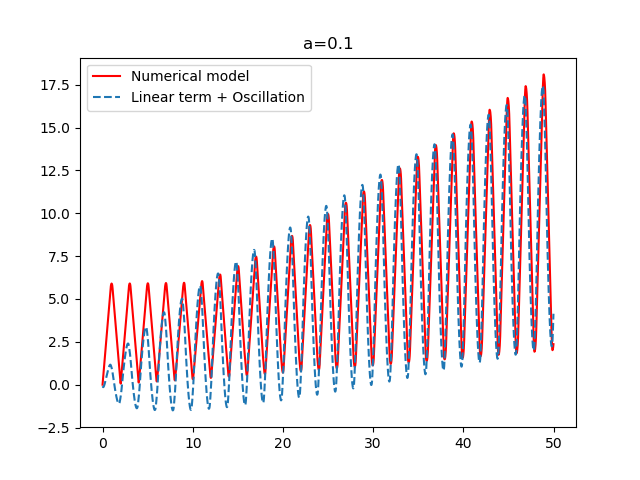}
\caption{Comparison between $|S_t|$ and $2\arcsin(a) t + t\sum _{k=0}^N J_k^{\rm approx}(t)$ for $N=10$, $t\in[0,50]$
 and $a \in \{0.1,0.3,0.7,0.9\}$.}
\end{figure}

First, we can observe that $\vert S_t \vert$ admits a linear asymptotic as $t$ grows to $+\infty$. Then, the oscillations are of period 2 with a phase independent of $a$ (see Figure \ref{frontdonde2}). One may remark the following points.
\begin{enumerate}
\item For $a=0$, the family of curves $(S_t)_t$ are concentric circles and $\vert S_t \vert$ is of period 2.
\item For $a=0.5$, the terms $J_k^{\rm approx}(t)$ vanish for any $t$ and, in this case, this expansion is not able to capture the oscillating part of $t\mapsto \vert S_t\vert$. 
\item The terms $\vert J_k^{\rm approx}(t)\vert$ are bounded by $Ck^{-5/2}t^{-1/2}$, where $C$ is a constant. For $t$ fixed, this ensures the (fast) convergence of the serie $\sum_k J_k^{\rm approx}(t)$ and then the amplitude of $\displaystyle t \mapsto t\sum _{k\in \N} J_k^{\rm approx}(t)$ is of order $t^{1/2}$ (see Figure \ref{frontdonde3}). 
\end{enumerate}

\section{A short proof of   the Arcsinus formula}

In the paper \cite{Vi-20}, the author was able to prove by elementary calculations the
\begin{theo}\label{theo:maindisk}
If $X$ is the unit disk, 
$|S_t|=\lambda (q) t+\mathrm{o}(t)$ as $t\to +\infty$ with radius $1$ with  $\lambda (q)=2 \arcsin a $ where  $a$
is the distance from $q$ to the center of the disk. 
\end{theo} 

We will reprove it using tools which will be extended to integrable geodesic flows in a forthcoming paper. 
For this, we will prove an integral formula:
\begin{theo}\label{theo:intdisk}
Let $\psi $ be the function periodic of period $1$ whose restriction to $[0,1]$ is given by $\psi (\theta )=|2\theta -1| $.
We have
\begin{equation}\label{equ:xi} |S_t|= t \Sigma _\pm  \int _{I_{\ga_0}} \psi \left(\theta _2^{\pm ,q} (\xi) -
 \frac{t}{2 \sin \xi }\right)d\xi + \mathrm{O}(1)\end{equation}
 as $t\to+\infty$, where $I_{\ga_0}:=[\pi/2 -\ga_0, \pi/2 + \ga_0 ]$,
$\alpha _0 =\arcsin a$ and 
\[  \theta_2^{\pm , q} (\xi)=\ha \pm \frac{\sqrt{a^2-\cos^2 \xi}}{2\sin \xi } \]
This integral can also be written as an integral over $\T$:
\begin{equation}\label{equ:alpha}
 |S_t|= t  \int _\T \psi \left(\ha - \frac{a\cos\alpha  +t}{2\sqrt{1-a^2 \sin^2 \alpha }} 
 \right)\frac{a \cos \alpha}{\sqrt{1-a^2 \sin^2 \alpha}}d\alpha  + \mathrm{O}(1)
\end{equation}
 as $t\to+\infty$.
\end{theo}
Let us show how Theorem \ref{theo:maindisk} follows from Theorem \ref{theo:intdisk}. 
We  consider an integral
\[ I(t)= \int _{I_{\ga_0}} \psi \left(\theta (\xi )-\frac{t}{2 \sin \xi }\right)d\xi \]
with $\theta $ smooth. 
We first approximate uniformly $\psi $ by a sequence of trigonometric polynomials
$\psi_N(u)=\sum _{|n|\leq N} a_n {\rm exp}(2i\pi n u ) $ with $a_0=\int_0^1 \psi (\theta ) d\theta =\ha $.
This way we get 
\[ I_N(t)=2 \ga _0 + \sum _{|n|\leq N,~n\ne 0}a_n   \int_{I_{\ga_0}}e^{2i\pi n\theta (\xi)}e^{-2i\pi nt/ \sin \xi } d\xi \]
It follows from the stationary phase approximations that all these integrals tend to $0$ as $t\ra \infty $, Theorem \ref{theo:maindisk} follows.

{\it Proof of Theorem \ref{theo:intdisk}.--}

We will first parametrize the dynamics using angle coordinates on tori. Let us denote by  $m(s)=(\cos s,\sin s) $ on the circle
and by $\vec{u}_s $ the vector $\overrightarrow{0m}(s)$.
Let us introduce a set of coordinates. In what follows, we parametrize 
 the 2D-submanifold   of  the phase space  consisting of oriented chords joining
a  point $m(s) $
to $ m(s+2\xi ) $ with speed $1$ by   $\xi \in ]0,\pi [$.
 Changing the orientation of the chords  moves $\xi $ into $\pi -\xi$. 
 For $\xi \in ]0,\pi [$ and $r\in [0, 2 \sin \xi ]$, we define
$F_\xi (s,r)= m(s) + r \vec{u}_{s+\xi +\pi/2}$.
This describes the chord $C_\xi $  between $m(s) $ and $ m(s+2\xi )$.
 The function $F_\xi$ is extended as a function on $\R^2$ periodic with respect to the lattice
$L_\xi $ spanned by the vectors $(2\pi,~ 0)$ and $(2\xi ,~ -2 \sin \xi )$. 
The function  $F_\xi $ is continuous, but only piecewise smooth. 
The pull-back under $F_\xi $ on $\R^2$ of the billiard dynamics is generated by the vector $\pa_r $. 

The coordinates $(s,r )$ range over a torus $\R^2/L_\xi $.  In order to continue the computation, we need to fix 
the lattice $\Z^2$. For that we introduce the linear map
$M_\xi:\R^2 _{\theta_1,\theta_2}\ra \R^2_{s,r} $ sending the canonical basis of $\Z^2$ onto the previous basis of $L_\xi $. 
The dynamics on the torus $\R^2/L_\xi $ is the image of $\pa_r $ under  $M_\xi ^{-1}$; let us denote it by $V$.
We get 
\[ V= \frac{1}{2\pi \sin \xi }(\xi \pa _{\theta_1} - \pi \pa _{\theta_2} )\]
Then, we need to compute the Euclidean norm of $F' _\xi (M_\xi (\pa _\xi V))$.
We have
\[ \pa _\xi V= \frac{-\cos \xi}{2\pi \sin ^2 \xi }(\xi \pa _{\theta_1} - \pi \pa _{\theta_2} )+ \frac{1}{2\pi \sin \xi}\pa_{\theta_1} \]
Hence
\[ M_\xi (\pa _\xi V)= \frac{1}{\sin \xi }\left(-\cos  \xi \pa _r + \pa _s \right) \]
Then 
\[ F'_\xi (\pa_r)= \vec{u}_{s+\xi +\pi/2},~  F'_\xi (\pa_s)=\vec{u}_{s+\pi/2}-r\vec{u}_{s+\xi} \] 
This gives
\[ \| \pa _\xi V   \|= \frac{|r-\sin \xi |}{\sin \xi}\] 
As could have been anticipated, this length vanishes on the caustic!
We now take the pull back of $  \|\pa _\xi  V  \|$ under $M_\xi$ and get $|2\theta_2-1|$.

Let us parametrize the chords starting from $q$ by the angle  $\ga \in \T$  defined
by $\ga:= \langle q,C_\xi \rangle$. We get $\cos\xi = a \sin \ga $. Hence $\xi $ is the smooth function
$\xi(\ga)=\arccos (a\sin\ga) $. The length $\vert S_t\vert$ is given by
$$
\vert S_t = \int_{\mathbb{T}}\Vert \frac{d}{d\alpha}(\phi_t(\vec{u}_\alpha)\Vert d\alpha
$$
where $\phi_t$ is the geodesic flow. Let us denote by $\theta (\ga)$ the coordinates of $q$ in $\T^2_\theta $. We get, using the parametrization of the flow on the tori $\mathbb{T}_\theta$,
\begin{eqnarray*}
\vert S_t \vert &=&
\int_{\mathbb{T}}\Vert(F'_{\xi} \circ M_{\xi})_{\theta(\alpha)+tV(\alpha)}(\theta'(\alpha)+tV'(\alpha))\Vert d\alpha,\\
&=&
t \int_{\mathbb{T}}\Vert(F'_{\xi} \circ M_{\xi})_{\theta(\alpha)+tV(\alpha)}(V'(\alpha))\Vert d\alpha + \mathrm{O}(1)
\end{eqnarray*}
as $t\to +\infty$. We rewrite the integral in terms of $\xi $, using $\cos \xi = a\sin \alpha $ and 
$\theta _2 (\xi)=\ha \pm  \frac{\sqrt{a^2 -\cos^2 \xi}}{2 \sin \xi } $
with $+$ if $\ga \in [\pi/2,3\pi/2]$ and $-$ otherwise. From this follows the result.

\section{Local asymptotics of the length}
In this section, we describe the asymptotics of the length of the intersection of the wave front with a smooth domain $K$
included in the disk $\D $. 
We have
\begin{theo}
We have
\[ l(S_t\cap K)  \sim \frac{2t}{\pi}\int_K \Psi\left(\sqrt{x^2+y^2}\right) |dx dy|  \]
as $t\to+\infty$, where
\[ \Psi(r)=  \frac{\min (r,a)}{\sqrt{1 -\min (r,a)^2}}\]

\end{theo}
Note that the function $\Psi$ is continuous, vanishes at $r=0$ and is constant for $a\leq r \leq 1$.
This implies that the density of the wave front is smaller near the center of the disk.

{\it Proof.--}
Let $\phi \in C( \D,\R^+)$, we want to calculate
the asymptotics of the length $|S_{t,\phi} |$ of $S_t$ computed in the metric $\phi^2 {\rm Eucl }$.
Following the proof of Theorem \ref{theo:maindisk},
we get $ |S_{t,\phi} |/t \ra \lambda (q,\phi ) $ as $t\ra +\infty $, with
\[ \lambda (q,\phi ) =2\int _{I_{\ga_0}}\int_{\T^2}|2\theta _2-1| \phi\circ G (\theta, \xi) |d\xi d\theta| \]
with $G(\theta ,\xi)=F_\xi \circ M_\xi (\theta )$.
We will first make the change of variable $(\theta,\xi)\ra (s,r,\xi)$ whose Jacobian is
$4\pi \sin \xi $. This gives
\[ \lambda (q,\phi ) =\frac{1}{2\pi}\int _{I_{\ga_0}}\int_{\R^2/L_\xi}\left|\frac{r-\sin \xi}{\sin^2 \xi}\right|
 \phi\circ F_\xi  (s,r) |d\xi ds dr | \]
Finally, we pass from $(s,r) $ to $(x,y)$. We have 
$|dxdy |= |r-\sin \xi | |ds dr |$. The domain of integration is $\rho=\sqrt{x^2+y^2}\geq \cos \xi $
which is covered twice by the torus ${\R^2/L_\xi}$, we get hence
\[ \lambda (q,\phi ) =\frac{1}{\pi}\int _{I_{\ga_0}}\int_{\cos \xi \leq \rho}\left|\frac{1}{\sin^2 \xi}\right|
 \phi (x,y) |d\xi dxdy  | \]
An elementary calculus gives then
\[ \lambda (q,\phi ) =\frac{2}{\pi}\int_\D \Psi (\rho) \phi (x,y) |dx dy |\]
The result follows then  by approximating the characteristic function of $K$ by continuous fonctions.
\hfill$\square $

\section{Oscillations of the length}

The numerical computations of the second author in \cite{Vi-20} show clearly some regular oscillations of the length $|S_t|$
around the linear asymptotics. These oscillations are given in the

\begin{theo}\label{theo:osc}  The following expansion holds:
\[ |S_t|= 2\ga _0 t + t\sum _{n=0}^\infty  J_n^{\rm approx}(t)+\mathrm{O}(1) \]
as $t\to +\infty$, with
\[ J_n^{\rm approx}(t)= \frac{-8\sqrt{2}}{\pi^2 (2n+1)^{5/2}\sqrt{t}}\cos \left( (2n+1)\pi a \right) \cos \left( 
\pi \left( (2n+1)t +\frac{1}{4}\right) \right)\]
The oscillations have an amplitude of the order of  $\sqrt{t}$,
are periodic of period $2$.
If $a=\ha $, we get $|S_t|= \pi t/3 + \mathrm{O}(1) $ as $t\to +\infty$. 

\end{theo}
We start from the formula given by Equation (\ref{equ:xi}):
\[ |S_t|= t \Sigma _\pm  \int _{I_{\ga_0}} \psi \left(\theta _2^{\pm ,q} (\xi) - \frac{t}{2 \sin \xi }\right)d\xi + \mathrm{O}(1)\]
as $t\to +\infty$, where $I_{\ga_0}:=[\pi/2 -\ga_0, \pi/2 + \ga_0 ]$ and  $\psi $ restricted to $[0,1]$ is given by $\psi (\theta )=|2\theta -1| $
 and $\psi $ is
periodic of period $1$.
We have
\[  \theta_2^{\pm , q} (\xi)=\ha \pm \frac{\sqrt{a^2-\cos^2 \xi}}{2\sin \xi } \]

The idea is to start with the Fourier expansion of $\psi $ and then to apply the stationary phase asymptotics.

We have
\[ \psi (\theta )=\ha  + \sum _{n\in \Z }\frac{2}{\pi^2 (2n+1)^2}e^{2(2n+1) i\pi \theta}\]
We  need to evaluate the integrals
\[ I_n (t)=\frac{-4}{((2n+1)\pi)^2}
\int _{I_{\ga_0}} \cos \left( (2n+1) \pi \frac{\sqrt{a^2-\cos^2\xi}}{\sin \xi }\right)e^{-i\pi (2n+1) \frac{t}{ \sin \xi }} d\xi \]
and then we have
\[ |S_t|= 2\ga _0 t + t\sum _{n\in \Z} I_n(t)+\mathrm{O}(1) \]
as $t\to +\infty$. Note first  that the function $\cos \left( (2n+1)\frac{\sqrt{a^2-\cos^2 \xi}}{\sin \xi }\right)$ is smooth on $I_{\alpha _0}$  with 
a non vanishing derivative  at the boundaries. 
The non vanishing contributions come from the critical point $\xi =\pi/2 $
and the boundaries of $I_{\ga_0}$. The boundary contributions are $\mathrm{O}(1/t)$. They contribute to the $\mathrm{O}(1)$ remainder.
The  contribution of the critical point can be calculated using the formula (\ref{equ:stat}). 
We get an asymptotic for $J_n= I_n + I_{-n-1},~n=0,1, \cdots$ given by
\[ J_n (t)\sim J_n^{\rm approx}= \frac{-8\sqrt{2}}{\pi^2 (2n+1)^{5/2}\sqrt{t}}\cos \left( (2n+1)\pi a \right) \cos \left( 
\pi \left( (2n+1)t +\frac{1}{4}\right) \right)\]

The previous calculation is only formal. We need to control the remainder terms in a uniform way with respect to $n$.
Let us rewrite the integral $I_n$ as combination of integrals of the form
\[ \int _{I_{\ga _0}}e^{-i\pi(2n+1)t \left(\frac{1}{\sin \xi}-\frac{1}{t}\frac{\sqrt{a^2-\cos^2 \xi}}{\sin \xi }\right)} d\xi \]
and apply the stationary phase with the phase functions depending on $t$:
$\Phi_t(\xi)=\frac{1}{\sin \xi}-\frac{1}{t}\frac{\sqrt{a^2-\cos^2 \xi}}{\sin \xi }$.
This phase function is non degenerate and converges in $C^\infty $ topology to $\frac{1}{\sin \xi}$ as $t\ra \infty $.
Hence the remainder is $\mathrm{O}\left( (nt)^{-3/2}\right)$ as $t\to +\infty$, uniformly with respect to $n$.

\begin{appendix}

\section{Stationary phase}
For this section, we refer the reader to \cite{GS-77}, chap. 1. 

We want to evaluate the asymptotics as $t\ra +\infty $ of integrals of the
form
\[ I(t):= \int _\T  e^{it S(x)} a(x) dx \]
where $S$ is a real valued smooth function.
 We assume that the critical points of $S$, ie the zeroes of $S'$, are
non degenerate, ie $S''(x)\ne 0$.
We will first assume that  $a\in C_o^\infty (\R)$ with only one critical point $x=0$ in the support of $a$. 
Then $I(t)$ admits a full asymptotic expansion given by 
\begin{equation} \label{equ:stat} I(t) =\frac{\sqrt{2\pi}e^{i\ge \pi /4}}{|t S'' (0)|^{\ha}}e^{itS(0)}\left(a(0)+ \mathrm{O}(t) \right)
  \end{equation}
as $t\to +\infty$, with $\ge = \pm 1$ depending on the sign of $S''(0)$.  
We will need some uniform estimates in the remainder term. This is provided by the following
\begin{prop}\label{prop:stat}  Let us consider the integrals
\[ I (t;S,a):= \int _\T  e^{it S (x)} a (x) dx \]
Let   $S_0$ be  a smooth real valued Morse function and $a_0$ be a  smooth function.
Let $S_\lambda $ and $a_\lambda $ be smoothly dependent of a real parameter $\lambda $.
Then, for $\lambda $ small enough, 
\[ I (t;S_\gl,a_\gl):=I_{\rm asympt}(t,\gl) + \mathrm{O}\left(t^{-3/2}\right)\]
as $t\to +\infty$, where the $\mathrm{O}$ is uniform and $I_{\rm asympt}(t,\gl)$ is the sum of terms given by the formula (\ref{equ:stat}) for all critical points
of $S_\gl$. 
\end{prop}
If $\gl $ is small enough, $S_\gl$ is still a Morse function. We localize the integrals near the critical points and apply 
the Morse Lemma with parameters. We are then reduced locally to the case where $S_\gl (x)=\pm x^2$. 
We apply then any proof of the stationary
phase approximation.

It will also be useful to consider the case of an integral on a closed interval $[c,d]$ with $c<d$.
\[ I(t):= \int _c^d  e^{it S(x)} a(x) dx \]
 Assuming that $S'$ does not vanish on the support of $a$  and that $a$ is $C^1$, we have 
\begin{equation} \label{equ:boundary} I(t)=\frac{1}{it}\left( \frac{a(d)e^{itS(d)}}{S'(d)} - \frac{a(c)e^{itS(c)}}{S'(c)}+O(t) \right) 
 \end{equation}
 as $t\to +\infty$.

Note that in both asymptotic formulae, the remainders ``$O(t)$'' are uniform if $S' $ (resp. $a'$) 
 is close to $S$ (resp. close to $a$) in $C^2$ topology.

The previous asymptotics extend to higher dimensional integrals.

\end{appendix} 
\bibliographystyle{plain}

\end{document}